\documentclass[11pt]{article}
\usepackage[T1]{fontenc}
\usepackage{graphicx}
\usepackage{amsmath}
\usepackage{amssymb}
\usepackage{algorithm}
\usepackage[noend]{algpseudocode}
\usepackage{authblk}
\usepackage{fullpage}

\begin{document}
\title{Efficient Algorithms for Energy-Aware Single-Machine Scheduling with Battery Storage}
\author[1,2]{\normalsize Luca Forte \thanks{luca.forte@ugent.be}}
\author[1,2]{\normalsize Stijn De Vuyst \thanks{stijn.devuyst@ugent.be}}
\author[1,2]{\normalsize Pieter Leyman \thanks{pieter.leyman@ugent.be}}
\affil[1]{Actemium Chair - Sustainable Energy, Department of Industrial Systems Engineering and Product Design, Ghent University, Ghent, Belgium}
\affil[2]{Industrial Systems Engineering (ISyE), Flanders Make, Kortrijk, Belgium}

%
%
%
%
%
\maketitle              
%


\begin{abstract}Energy-aware scheduling has become a central challenge in modern manufacturing environments. As industries increasingly aim to reduce operational costs and carbon emissions, aligning production activities with electricity tariffs is increasingly important. The integration of battery energy storage system (BESS) further enhances the potential to reduce energy purchase costs; however, incorporating such systems introduces complex interdependencies between job scheduling and BESS management decisions, significantly complicating the problem structure.\\In this work, we study the problem of minimizing the energy cost of executing a set of jobs on a single machine within a fixed time horizon, where electricity prices follow a Time-of-Use (TOU) tariff. In addition, we consider the presence of a BESS that can be charged from the grid and discharged during high-price periods to reduce overall energy costs. To efficiently address this problem, we propose and analyze two matheuristic algorithm variants designed to effectively coordinate production scheduling and BESS usage decisions.
\end{abstract}
%
%
%
\section{Introduction}\label{introduction_literature_review}
With the ongoing energy transition, the coordination of production scheduling and energy management decisions is becoming increasingly critical, as this will combine both a reduction of industrial operational costs and an increased support for the broad transition toward low-carbon energy systems. Energy-aware scheduling has therefore emerged as a central challenge in production scheduling: in contrast to classical machine scheduling problems, where objectives typically focus on makespan, energy-aware models explicitly account for time-varying electricity prices and energy consumption profiles. Under Time-of-Use (TOU) tariffs, electricity prices vary across predefined time intervals, creating incentives to shift energy-intensive activities toward low-price periods. Additionally, many facilities can also be equipped with a battery energy storage system (BESS). These systems enable electricity to be purchased and stored during low-price periods and used during high-price periods, thereby providing further opportunities for cost savings. From a modeling perspective, a BESS introduce a new layer of decisions: we need to decide when to charge, when to discharge, and how to allocate stored energy across jobs over time. As a result, production scheduling and BESS management decisions become strongly interdependent, which significantly enlarges the decision space and introduces additional combinatorial complexity beyond that of classical energy-aware scheduling.\\Several studies have proposed integrated models to coordinate machine operations with BESS decisions; however, important limitations remain. First, some works assume simplified electricity pricing structures. In particular, discretized TOU schemes with only three tariff levels (e.g., peak, mid-peak, off-peak) are typically employed, instead of real day-ahead price data \cite{Moon,Ghorbanzadeh}. Second, unrealistically simple assumptions are often made in BESS formulations: charging and discharging decisions are determined using greedy rule-based approaches or modeled without detailed operational constraints, such as charging/discharging limits and efficiency losses \cite{Dong,Ghorbanzadeh,Liu}. Third, some existing approaches rely exclusively on exact mathematical programming formulations. While these models provide valuable structural insights, scalability issues are rarely considered \cite{Karimi,Liu}. In order to address these gaps, this work studies the problem of minimizing the energy cost of executing a set of jobs on a single machine within a fixed time horizon under TOU tariffs, while accounting for sequence-dependent setup times and the presence of a BESS. To the best of our knowledge, this specific setting has not been explicitly investigated with the primary objective of cost minimization. The problem is formulated as a Mixed-Integer Linear Program (MILP) and two matheuristic algorithms are developed to provide high quality solutions within reasonable computational time.
\section{Problem Definition} \label{problem_definition}
We consider a single-machine scheduling problem with time-dependent electricity prices and a BESS. Let $N$ denote a set of jobs, indexed by $i \in \{1,...,|N|\}$, with 0 and $|N|+1$ denoting start and end dummy jobs respectively. Each job $i$ has a processing time $d_i$, measured in discrete time periods. During period $k \in \{1,...,d_i\}$, job $i$ requires $E_{i, k}$ units of energy to be processed. The electrical power needed to process each operation is sourced from the public grid at day-ahead market prices. These prices are publicly available and are assumed to be known in advance for the entire planning horizon. Moreover, it is assumed that energy acquisition does not have any availability constraints. All jobs must be scheduled within a time horizon consisting of $T$ 15-minutes long intervals; the electricity price during period $t$ is denoted by $P_t$, with $t \in \{1,...,T\}$. If job $i$ is scheduled to start at time $t$, then its total energy cost equals $C_{i, t} = \sum_{k=1}^{d_i}P_{t+k-1}E_{i, k}$. Finally, for each pair of jobs $(i,j)$, with $i\neq j$, a sequence-dependent setup time $v_{i,j}$ is introduced; it is assumed that the machine remains idle during this setup and does not consume any energy. Furthermore, both job preemption and parallel processing are not allowed.\\
In addition to the operational scheduling environment described above, storage capabilities are incorporated through a BESS with the following features: a maximum capacity level $B^{cap}$, a maximum charging rate $B^{c, max}$, a maximum discharging rate $B^{d, max}$, charging and discharging efficiencies $\eta^c$ and $\eta^d$. Importantly, it is assumed that the BESS is not allowed to sell electricity to the grid and to perform simultaneous charging and discharging operations during a certain interval $t$.\\A mathematical formulation of the problem, called \textsc{milp-solution}, is provided below, along with an overview of the variables and parameters not discussed above (cf.\ Table \ref{tab1}). The objective function \eqref{obj1} represents the total electricity cost of the schedule, and the model seeks to minimize it. The first term is the total electricity cost of the schedule, increased with an additional cost when electricity is purchased for storage in the BESS, and decreased based on the amount discharged from the BESS. 
Constraints \eqref{c1} and \eqref{c2} ensure that each job has exactly one successor and one predecessor. Constraints \eqref{c3}--\eqref{c5} define the variable $s_i$ and employ it to enforce the required setup times between consecutive jobs. Constraints \eqref{c6} guarantee that each job is scheduled exactly once. Constraints \eqref{c7}--\eqref{c9} model the BESS dynamics, including charging and discharging decisions and state-of-charge evolution over time. Finally, constraints \eqref{c10} impose an additional bound on the discharging capability of the BESS relative to the active load: at time $t$, if no job is being processed, discharging the BESS would not provide any revenue and result in a waste of energy, so the discharge is bounded by 0. For the same reason, if a job is being processed, the discharge is bounded by the energy requirement of the job period executed during $t$. Constraints \eqref{c4} and \eqref{c9} rely on Big-M formulations; a constant value of $M = 1000$ was employed.\\As shown in \cite{Fang}, the problem of minimizing the cost of scheduling a set of $|N|$ jobs with a duration $d_i$ and a constant power demand $q_i$ on a single machine under electricity tariffs is strongly NP-hard. The problem addressed in this paper generalizes that setting in several important directions and is therefore significantly more complex. First, jobs are characterized by time-dependent power consumption rather than a constant energy demand. Second, sequence-dependent setup times are introduced between consecutive jobs. Third, an on-site BESS is integrated into the system, requiring joint optimization of production scheduling and energy storage decisions. Therefore, a feasible solution is defined by a permutation $\pi$ of the job set $N$, a vector of start times $s \in \mathbb{N}_+^N$ and a BESS operations schedule over the time horizon. Preliminary experiments indicate that solving \textsc{milp-solution} to optimality is computationally intractable for real-world instances; therefore, two metaheuristic approaches are proposed to efficiently obtain high-quality solutions.
\begin{table}[h!]
\caption{Summary of symbols and their definitions.}\label{tab1}
\setlength{\tabcolsep}{15pt}
\renewcommand{\arraystretch}{1.3}
\begin{tabular}{l p{7.5cm}}
\hline
Parameters & Definition \\
\hline
$v_{i,j}$ & Setup time from job $i$ to job $j$\\
\hline
Decision variables & Definition\\
\hline
$b_t^c$ & Amount of energy purchased to charge the BESS during $t$\\
$b_t^d$ & Amount of energy discharged from the BESS during $t$\\
$x_{i,t}$ & 1 if job $i$ is scheduled during $t$, 0 otherwise\\
$z_{i,j}$ & 1 if job $j$ follows job $i$, 0 otherwise\\
$s_i$ & Start time of job $i$\\
$y_t$ & 1 if the BESS is charging during $t$, 0 if the BESS is discharging during $t$ \\
\end{tabular}
\end{table}

\begin{align}
\text{min} & \sum_{i=1}^{|N|} \sum_{t=1}^{T}C_{i,t} x_{i,t} + \sum_{t=1}^{T}P_t b_t^c - P_t b_t^d\label{obj1} \\
\text{s.t.} \quad
& \sum_{j=1, j\neq i}^{|N|+1} z_{i,j} = 1, \quad \forall i = 0,...,|N|\label{c1} \\
& \sum_{j=0, j\neq i}^{|N|} z_{j,i} = 1, \quad \forall i = 1,...,|N|+1\label{c2} \\
& s_i=\sum_{t=1}^{T-d_i} tx_{i,t}, \quad \forall i = 1,...,|N| \label{c3} \\
& s_i + d_i + v_{i, j} - M(1-z_{i,j}) \leq s_j, \quad \forall i=0,...,|N|, \forall j=1,...,|N|+1, i \neq j \label{c4} \\
& 0 \leq s_i \leq T-d_i, \quad \forall i = 1,...,|N| \label{c5} \\
& \sum_{t=1}^{T-d_i}x_{i,t}=1, \quad \forall i = 1,...,|N| \label{c6} \\
& b_t=b_{t-1} + \eta^cb_{t-1}^c -\frac{b_{t-1}^d}{\eta^d}, \quad \forall t=1,...,T \label{c7} \\
& b_t \leq B^{cap}, b_t^c \leq B^{c, max}, b_t^d \leq B^{d, max}, \quad \forall t=1,...,T \label{c8} \\
& b_t^c \leq M y_t, b_t^d \leq M(1-y_t), \quad \forall t=1,...,T \label{c9} \\
& b_t^d \leq \sum_{i=1}^{|N|} \sum_{s=max(1, t-d_i)}^{t} x_{i,s}E_{i, t-s}, \quad \forall t=1,...,T \label{c10}\\
& x_{i, t} \in \{0, 1\}, \quad \forall i=1,...,|N|, \forall t=1,...,T-d_i\label{c11} \\
& z_{i, j}\in \{0, 1\}, \quad \forall i=0,...,|N|, \forall j=1,...,|N|+1, i \neq j\label{c12} \\
& y_t\in \{0, 1\}, \quad \forall t=1,...,T\label{c13} \\
& b_t^c, b_t^d, b_t \geq 0, \quad \forall t=1,...,T\label{c14}\\
& s_0 = 0, s_{|N|+1} = T, b_0 = 0\label{c15}
\end{align}

\section{Methods}\label{methods}
\subsection{MH-SEQ-MILP Algorithm}\label{mh_seq_milp}
The first approach is a matheuristic combining Iterated Local Search (ILS) with a reduced MILP model. The ILS explores the space of job permutations $\Pi$, i.e., all possible sequences of the $|N|$ jobs. Then, for a given sequence $\pi \in \Pi$, a modified MILP model, denoted \textsc{milp-schedule}, is solved to jointly determine the optimal job start times and BESS operations. Since the job order is fixed by $\pi$, \textsc{milp-schedule} is a reduced and computationally more efficient version of \textsc{milp-solution}, as all the $z_{i, j}$ variables and respective constraints can be dropped. The overall algorithm described above is referred to as \textsc{mh-seq-milp}. For the ILS, a solution is represented as a permutation vector $\pi = (\pi(1),...,\pi(|N|))$, where $\pi(k)$ denotes the job processed in position $k$. Idle times between jobs are not explicitly represented, as they are determined by \textsc{milp-schedule} during optimization. The initial solution is generated randomly. The perturbation operator works on consecutive blocks of jobs within the sequence. Given a fixed block size $b$, the sequence is divided into contiguous blocks consisting of $b$ jobs each; for each block, the total cost is computed as the sum of the individual job costs. The block with the highest cost is then selected, and the jobs within it are randomly shuffled to generate a perturbed sequence; based on preliminary testing, a block size of $\lfloor{\frac{|N|}{5}}\rfloor$ was chosen for each instance. Finally, the acceptance criterion follows a restart-based logic: the search restarts from a random sequence after 5 iterations without improvement. While the values for both the block size and the number of iterations parameters were fixed to obtain preliminary results, they are candidates for further optimization: we intend to refine these values and explore adaptive tuning mechanisms.

\subsection{MH-SEQ-HYBRID Algorithm}\label{mh_seq_hybrid}
The second approach, denoted \textsc{mh-seq-hybrid}, decomposes the problem into three phases. As in \textsc{mh-seq-milp}, an ILS procedure first determines a job sequence $\pi$. Given this sequence, a backward Dynamic Programming (DP) algorithm, \textsc{dp-timing}, assigns job start times by minimizing cumulative electricity costs, thereby fixing the timing of the  electricity demand over the horizon $T$. Algorithm \ref{dp} provides the pseudocode for \textsc{dp-timing}. The entry $opt[i][t]$ represents the minimum cost for scheduling the subsequence of jobs $\{i,...,|N|\}$ starting from time $t$; the function $C(i, t)$ computes and returns the cost of scheduling job $i$ at time $t$ according to the definition of $C_{i,t}$ in Section 2; if the job's execution exceeds $T$, $C(i,t)$ returns $\infty$. Subsequently, another reduced mathematical model, \textsc{milp-bess}, optimizes BESS charging and discharging decisions for the fixed schedule. Since job start times are already determined, the resulting model is significantly more efficient than \textsc{milp-schedule}. Unfortunately, this decomposition introduces a structural limitation. Because job timing and BESS optimization are solved sequentially rather than jointly, there is no theoretical guarantee that, for a given sequence $\pi$, the combination of \textsc{dp-timing} and \textsc{milp-bess} yields an optimal solution. Therefore, this method trades optimality guarantees for computational efficiency, a compromise that will be experimentally evaluated.

\begin{algorithm}
\caption{\textsc{dp-timing}}\label{dp}
\begin{algorithmic}[1]
    \Require $|N|, T, \text{durations}\ d, \text{setup times}\ v$
\State $opt[i][t] \gets \infty \quad \forall i \leq |N|, t \leq T$
\For{$i\gets |N|$ to $1$}
\For{$t\gets T$ to $1$}
\If{$i = N$}
\State $t_{\text{next}}\gets t+d_i$
\Else
\State $t_{\text{next}}\gets t+d_i+v_{i, i+1}$
\EndIf
\State $opt[i][t] \gets \min(\text{COMP\_V}(i, t+1, opt), \text{COMP\_V}(i+1, t_{\text{next}}, opt])+C(i, t))$
\EndFor
\EndFor
\State \Return opt
\end{algorithmic}
\end{algorithm}

\begin{algorithm}
\caption{\textsc{comp\_v}}\label{p}
\begin{algorithmic}[1]
\Require $i, t, opt$
\If{$i \geq |N|$}
\State \Return $0$
\ElsIf{$t \geq T$}
\State \Return $\infty$
\Else
\State \Return $opt[i][t]$
\EndIf
\end{algorithmic}
\end{algorithm}

\section{Results}
\subsection{Dataset}The test instances employed in the experiments are derived from \cite{Mencaroni} and were modified to ensure feasibility after introducing setup times. Based on these adapted instances, four groups of five instances were generated, differing in the number of jobs, planning horizon length and slack time characteristics. An overview of the instance groups is reported in Table \ref{tab2}.
\begin{table}[t]
\centering
\caption{Computational results for different instance groups. 1 d = 1 day, 3 d = 3 days, 6 d-l-s = 6 days with low slack, 6 d-h-s = 6 days with high slack. For the \textsc{milp-solution} cost row, brackets report the average MIP gap; for the \textsc{mh-seq-milp} and \textsc{mh-seq-hybrid} cost rows, brackets report the average number of iterations.}\label{tab2}
\setlength{\tabcolsep}{6.5pt}
\renewcommand{\arraystretch}{1.3}
\begin{tabular}{l|l|l|l|l}
 & 1 d & 3 d & 6 d-l-s & 6 d-h-s\\
\hline
Average $|N|$ &10&30.6&63.2&31.9\\
\hline
$T$ &96&288&576&576\\
\hline
\textsc{milp-solution} cost & $\mathbf{3.62\times10^6}$ & $\mathbf{11.82\times10^6}$ & $\mathbf{36.22\times10^6}$ & $11.09\times10^6$\\
&(34.39\%)&(175.14\%)&(209.22\%)&(2763.83\%)\\
\hline
\textsc{mh-seq-milp} cost &  $3.96\times10^6 $&$13.19\times10^6$&$38.61\times10^6$&$11.89\times10^6$\\
&(678.4)&(90.6)&(11)&(0.6)\\
\hline
\textsc{mh-seq-hybrid} cost &$4.14 \times 10^6$& $12.52\times10^6$&$36.39\times10^6$&$\mathbf{10.53\times10^6}$\\
&(2469.6)&(795.6)&(380.6)&(458.2)\\
\end{tabular}
\end{table}
\subsection{Experimental design} 
To obtain a benchmark reference, \textsc{milp-solution} was executed on each instance with a time limit of one hour. The two metaheuristic approaches, \textsc{mh-seq-milp} and \textsc{mh-seq-hybrid}, were executed with a time limit of 60 seconds per instance. This setup allows us to evaluate whether the two metaheuristics can produce competitive solutions within short computational times and to directly compare their relative performance. Table \ref{tab2} reports the average results over the five instances of each group. The experiments were conducted on a machine equipped with an Intel(R) Core(TM) Ultra 5 236V (2.10 GHz) processor, using Python version 3.13.6 and gurobipy version 12.0.3.
\section{Conclusions}
The results in Table \ref{tab2} indicate that \textsc{mh-seq-hybrid} consistently outperforms \textsc{mh-seq-milp} in terms of solution quality for all instance groups except the smallest one (column 1 d). This difference can largely be attributed to the significantly greater number of permutations explored by \textsc{mh-seq-hybrid} within the assigned time limit. Overall, \textsc{milp-solution} achieves the best performance across nearly all instance groups, with the only exception being the largest and high-slack setting (column 6 d-h-s), where \textsc{mh-seq-hybrid} performs better. An additional insight, particularly relevant from a practical perspective, concerns scalability under high-slack conditions. In the 6--days high-slack group (column 6 d-h-s), the performance of \textsc{mh-seq-milp} deteriorates significantly, as it was often unable to complete even a single iteration. Interestingly, performance is worse in this high-slack group than in the 6--days low-slack group (column 6 d-l-s), despite the latter containing substantially more jobs. This suggests that \textsc{mh-seq-milp} can be suitable for tightly constrained instances but struggles when slack time increases.
\section*{Acknowledgments}

This study was funded by Ghent University's Special Research Fund (BOF) under grants BOFBAF2024097201, BOFBAF2025002101, and by the Actemium chair on Sustainable Energy.
%
%
%

\begin{thebibliography}{8}
\footnotesize
\bibitem{Dong}
Dong, J., Ye, C.: Green scheduling of distributed two-stage reentrant hybrid flow shop considering distributed energy resources and energy storage system. Computers \& Industrial Engineering \textbf{169}(108146) (2022)
Author, F.: Article title. Journal \textbf{2}(5), 99--110 (2016)

\bibitem{Fang}
Fang, K., Uhan, N. A., Zhao, F., Sutherland, J. W.: Scheduling on a single machine under time-of-use electricity tariffs. Annals of Operations Research \textbf{238}(1), 199--227 (2016)

\bibitem{Ghorbanzadeh}
Ghorbanzadeh, M., Ranjbar, M.: Energy-aware production scheduling in the flow shop environment under sequence-dependent setup times, group scheduling and renewable energy constraints. European Journal of Operational Research \textbf{307}(2), 519--537 (2023)

\bibitem{Karimi}
Karimi, S., Kwon, S.: Comparative analysis of the impact of energy‐aware scheduling, renewable energy generation, and battery energy storage on production scheduling. International Journal of Energy Research \textbf{45}(13), 18981--18998 (2021)

\bibitem{Liu}
Liu, C. H.: Mathematical programming formulations for single-machine scheduling problems while considering renewable energy uncertainty. International Journal of Production Research \textbf{54}(4), 1122--1133 (2016)

\bibitem{Mencaroni}
Mencaroni, A., Leyman, P., Raa, B., De Vuyst, S., Claeys, D: Towards net-zero manufacturing: Carbon-aware scheduling for GHG emissions reduction. Journal of Cleaner Production \textbf{529}(146787), (2025)

\bibitem{Moon}
Moon, J. Y., Park, J.: Smart production scheduling with time-dependent and machine-dependent electricity cost by considering distributed energy resources and energy storage. International Journal of Production Research \textbf{52}(13), 3922--3939 (2014)

\end{thebibliography}
%

\end{document}